\documentclass[11pt,a4paper,reqno]{amsart}
\usepackage[utf8]{inputenc}
\usepackage[T1]{fontenc}
\usepackage{amsmath}
\usepackage{amsfonts}
\usepackage{amssymb,amsthm}
\usepackage{graphicx}
\usepackage{mathrsfs}
\usepackage{tikz}
\usepackage{hyperref}
\usepackage{appendix}
\usepackage{enumitem}
\usepackage{color}
\usepackage{etoolbox}
\usepackage{comment}
\newtheorem{thm}{Theorem}[section]

\newtheorem{pro}[thm]{Proposition}
\newtheorem{cor}[thm]{Corollary}
\newtheorem{lem}[thm]{Lemma}
\newtheorem{rem}[thm]{Remark}
\newcommand{\R}{\mathbb{R}}
\newcommand{\N}{\mathbb{N}}

\def\qq#1{\qquad \mbox{#1}\quad}
\def\q#1{\quad \mbox{#1}\ }
\newcommand{\al}{\alpha}
\newcommand{\be}{\beta}
\newcommand{\De}{\Delta}
\newcommand{\de}{\delta}

\newcommand{\e}{\varepsilon}

\newcommand{\na}{\nabla}

\newcommand{\Om}{\Omega}
\newcommand{\Omb}{\overline{\Om}}
\newcommand{\p}{\partial}

\newcommand{\vf}{\varphi}

\makeatletter
\@addtoreset{equation}{section} 
\patchcmd{\maketitle}{\@fnsymbol}{\@alph}{}{}  
\makeatother

\title{Regularity and explicit $L^{\infty}$ estimates 
for a class of nonlinear elliptic systems
}

{\small
\author{Maya Chhetri}
\address{University of North Carolina at Greensboro, Greensboro, NC, USA}
\email{m\_chhetr@uncg.edu}
\author{Nsoki Mavinga}
\address{Swarthmore College, Swarthmore, PA, USA}
\email{nmaving1@swarthmore.edu}

\author{Rosa Pardo}
\address{Universidad Complutense de Madrid, Madrid, Spain}
\email{rpardo@ucm.es}
}

\date{\today}

\begin{document}
\maketitle
{ }\makeatletter{\renewcommand*{\@makefnmark}{}
}

{ }\makeatletter{\renewcommand*{\@makefnmark}{}
}

\begin{abstract}
We use De Giorgi-Nash-Moser iteration scheme to establish 
that weak solutions to a coupled system of 
elliptic equations with critical growth  
on the boundary are in $L^\infty(\Om)$. Moreover, we  provide an explicit   $L^\infty(\Om)$- estimate of weak solutions with subcritical growth on the boundary, in terms of powers of $H^1(\Omega)$-norms,   
by combining the elliptic regularity of weak solutions with Gagliardo--Nirenberg interpolation inequality. 
\end{abstract}

\section{Introduction}
We consider a coupled system of the form	
\begin{align}
\label{sys:u}
-\De u+u & =0 \q{in}\Om\,,
&\frac{\p u}{\p\eta}=f(x,v)\q{on}\p\Om\,;\\
-\Delta v +v &=  0 \q{in}\Om\,, 
&\frac{\p v}{\p\eta}=g(x,u)\q{on}\p\Om\,;
\label{sys:v}
\end{align}
where $\Om \subset \R^{N} (N \ge 3)$ is a bounded domain with $C^{0,1}$ (Lipschitz) boundary 
$\p \Om$, and  $\p/\p\eta
:=\eta(x)\cdot\nabla$ denotes the outer normal derivative on the boundary $\p\Om$. We assume that  $f,\ g:\p\Om\times \R \to \R$ are   Carathéodory functions, that is,  $f(\cdot,s),\ g(\cdot,s)$ are measurable for all $s\in\R$,  and $f(x,\cdot),\ g(x,\cdot)$  
are   continuous for a.e. $x \in \p\Om$. Moreover, there exist positive constants $b_1, b_2$ such that for a.e. $x \in \p\Om$
\begin{align}\label{f:g:growth}
|f(x,s)|\le b_1(1+|s|^{p_2}),\qquad |g(x,s)|\le b_2(1+|s|^{p_1}), \quad \forall s \in \R
\end{align}
for  $p_1, p_2 >1$ satisfying additional conditions. 

 \par Our goal is to establish an $L^\infty(\Om)$ {\it a priori estimate} for weak solutions of \eqref{sys:u}-\eqref{sys:v}.  The $L^{\infty}$ boundedness of weak solutions of elliptic equations are important for improving the regularity of weak solutions, 
 which is less understood than in the case of  Dirichlet boundary conditions. See \cite{Quittner_Souplet_2012}, where authors discussed challenges  associated to the parabolic problems involving nonlinearities in the interior vs. nonlinearity at the boundary. \\

 \par We say that a weak solution $(u,v)$ of \eqref{sys:u}-\eqref{sys:v} has $L^\infty(\Om)$ {\it a priori estimate} if  
 $u,v \in L^\infty(\Om)$ and $\|u\|_{L^\infty(\Om)}\leq M_1,\ \|v\|_{L^\infty(\Om)} \leq M_2$, where  $M_i=M_i(u,v,f,g,\Om)$ for $i=1,2$. 
\par
\begin{figure}[ht]
\begin{center}
\end{center}
\vspace{1cm}
$$
\begin{array}{ccc}
\kern -1cm N=3 
&\kern -1.3cm 3<N<4 
&\kern -1.3cm N\ge 4 \\
\kern -1cm
\includegraphics[width=5.5cm]{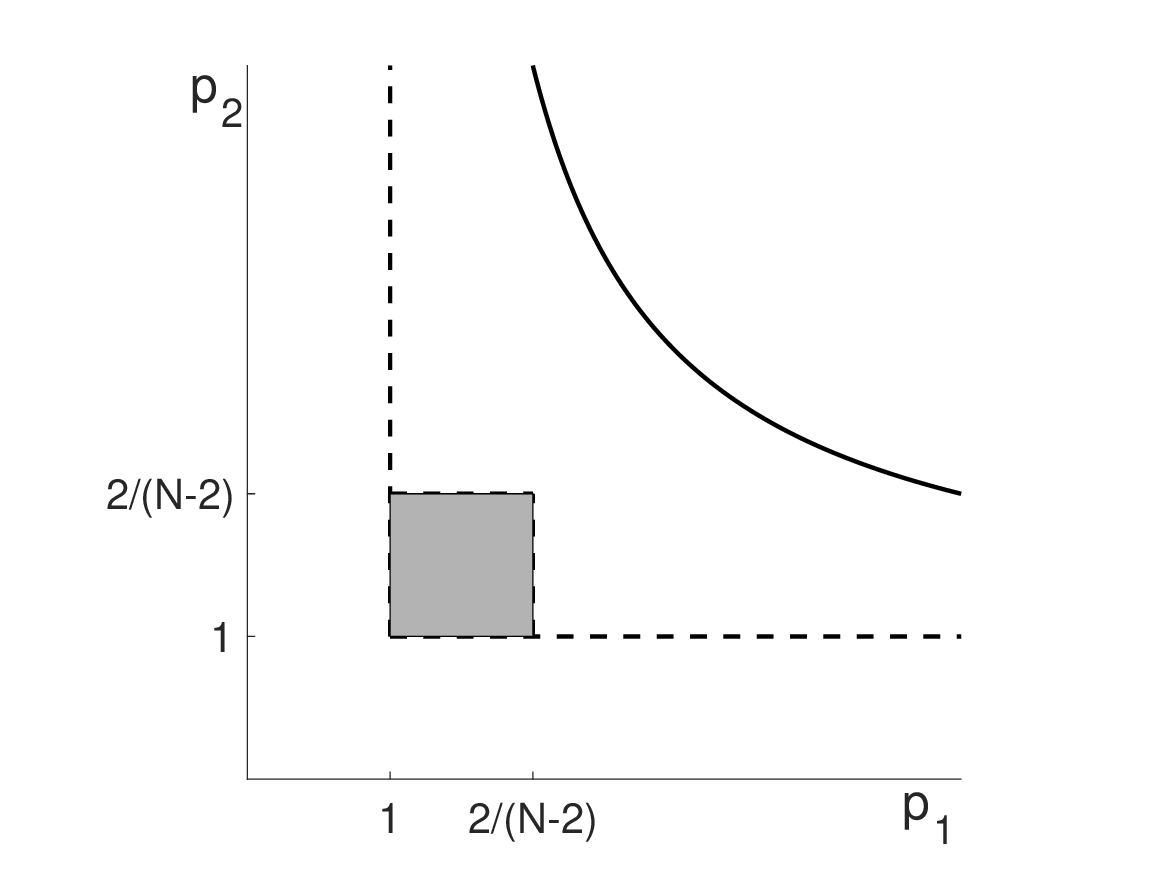} &
\kern -1.3cm
\includegraphics[width=5.5cm]{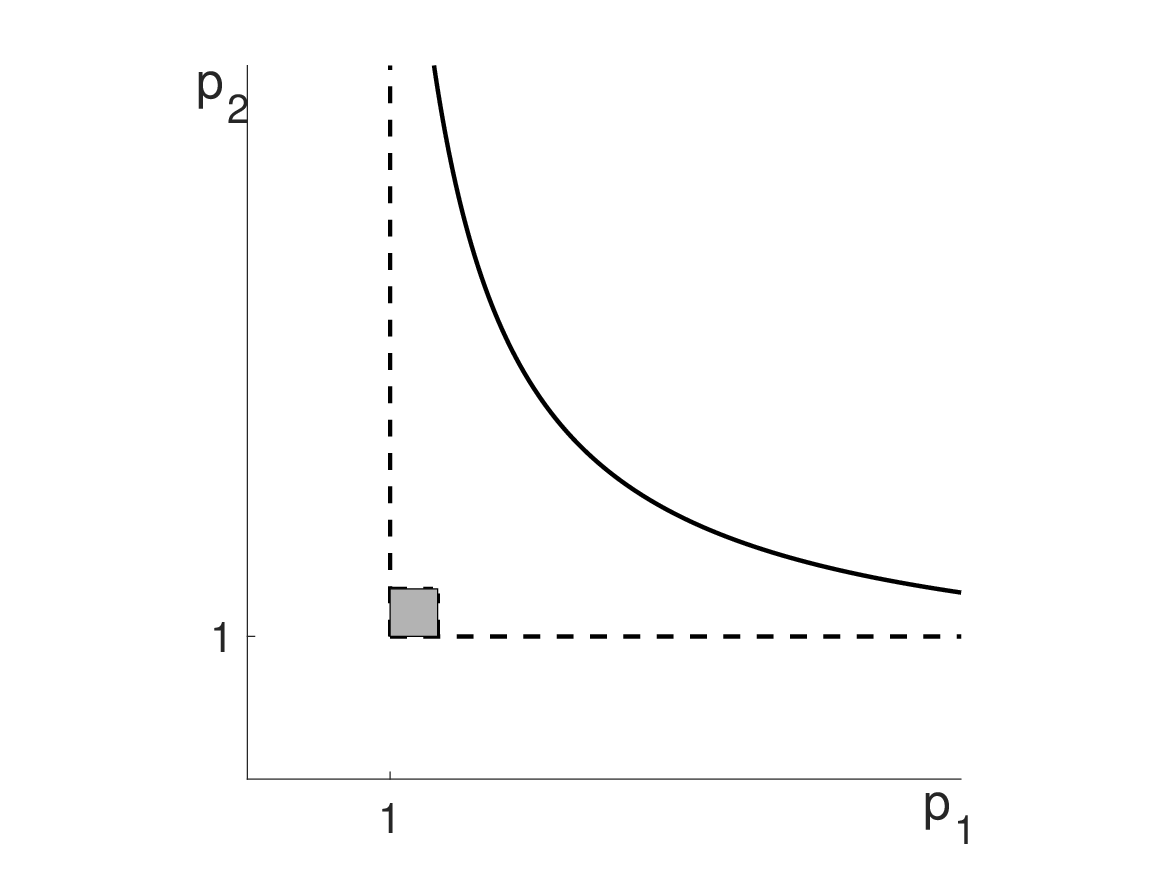} &
\kern -1.3cm
\includegraphics[width=5.5cm]{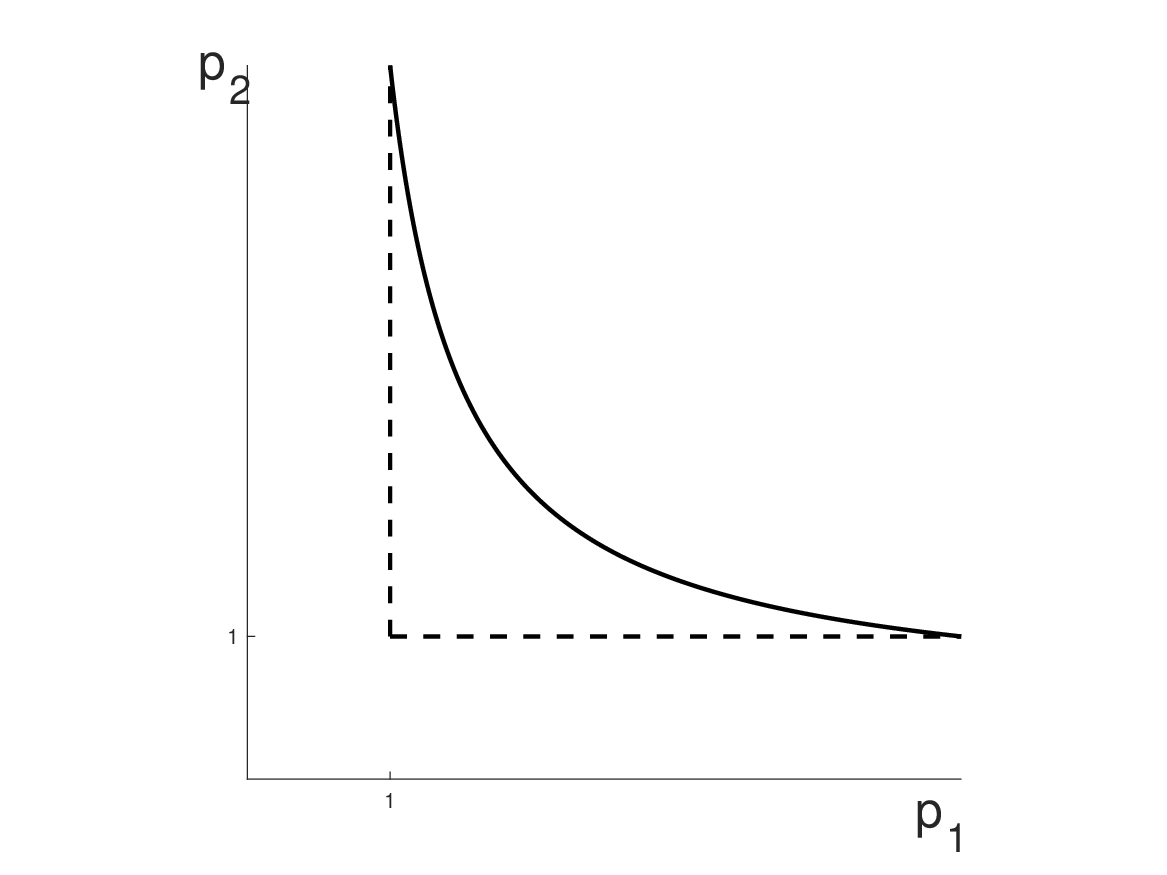}
\end{array}	
$$
\caption{Theorem~\ref{sys:th:reg} holds for $p_1,\ p_2>1$ under the hyperbola $
\frac{1}{p_1+1} + \frac{1}{p_2+1}= \frac{N-2}{N-1}
$  and the result of \cite{Marino-Winkert_systems_2020} holds on the shaded square $\big(1, \frac{N}{N-2}\big)\times \big(1, \frac{N}{N-2}\big)$, which is empty for any $N\ge 4$.  } 
 \label{fig:region}   
\end{figure}

 We first focus on proving the $L^{\infty}(\Om)$ regularity of weak solutions of \eqref{sys:u}-\eqref{sys:v} when $p_1, p_2$ satisfy 
\begin{equation}\label{sys:crit:hyp:2}
\frac{1}{p_1+1} + \frac{1}{p_2+1}\ge \frac{N-2}{N-1}; \qquad 1 < p_1\le p_2.
\end{equation}
Without loss of generality and by interchanging the roles of $p_1$ and $p_2$, the region given by \eqref{sys:crit:hyp:2} is under the hyperbola $\frac{1}{p_1+1} + \frac{1}{p_2+1}= \frac{N-2}{N-1}$, 
known as the {\it critical hyperbola}, see Figure~\ref{fig:region}. Such a  critical hyperbola is considered in \cite{Bonder-Pinasco-Rossi} for a Hamiltonian system. 

\par In the theorem below we establish the $L^{\infty}(\Om)$ regularity result for weak solutions to \eqref{sys:u}-\eqref{sys:v}. 
\begin{thm}[Regularity]
\label{sys:th:reg} 
Let $(u, v)$  be  a weak solution of  the system \eqref{sys:u}-\eqref{sys:v}. If the nonlinearities $f,\ g$ satisfy \eqref{f:g:growth} and \eqref{sys:crit:hyp:2} 
then  $u,v\in L^{q}(\p\Om)$ for any $1\le q<\infty$.
\smallskip

Additionally, $u, v \in C^\mu(\overline\Om)\cap W^{1,m}(\Om)$  for any $\mu \in (0, 1)$  and  $1<m<\infty$, with
\begin{equation}\label{infty:cont:u}
 \|u\|_{L^{\infty}(\p\Om)}\le \|u\|_{C(\Omb)}=\|u\|_{L^{\infty}(\Om)},   
\end{equation}
and
\begin{equation}\label{infty:cont:v}
 \|v\|_{L^{\infty}(\p\Om)}\le \|v\|_{C(\Omb)}=\|v\|_{L^{\infty}(\Om)}.    
\end{equation}
\end{thm}

\par Theorem~\ref{sys:th:reg} guarantees that if $(u,v)$ is a weak solution  to \eqref{sys:u}-\eqref{sys:v}, then
$u,v\in L^{\infty}(\Om)$. Next, to establish $L^\infty(\Om)$ {\it a priori estimate}, we provide explicit
expression for $M_i=M_i\left(\|u\|_{H^1(\Om)},\|v\|_{H^1(\Om)}\right)$,  $i=1,2$, 
which  is new for elliptic systems with nonlinearities on the boundary. 

\par 
For this, we consider  the region strictly below the critical hyperbola as follows:
\begin{equation}
\label{crit:hyp}
\frac{1}{p_1+1} + \frac{1}{p_2+1}> \frac{N-2}{N-1}; \qquad 1 < p_1\le p_2,
\end{equation}
see Figure~\ref{fig:region}. We prove the following result.
\begin{thm}
\label{sys:th:apriori:estim}
Let $f\,, g$  satisfy \eqref{f:g:growth}, \eqref{crit:hyp} and define 
\begin{equation}\label{de:0}
\de_0:=   
\frac{1}{p_1+1} + \frac{1}{p_2+1}-\frac{N-2}{N-1}>0 .
\end{equation}
If $(u, v)$ is a   weak solution of \eqref{sys:u}-\eqref{sys:v}, then there exist constants $C_0$ and $C_1$, independent of $u$ and $v$, such that 
\begin{equation}
\label{sys:est:u:L:inf}
\|u\|_{L^{\infty}(\Om)} \leq   C_0\, 
\bigg(1+\|u\|_{H^1(\Om)}^{A} \bigg)
\bigg(1+\|v\|_{H^1(\Om)}^{A} \bigg)
\end{equation}
and 
\begin{equation}
\label{sys:est:v:L:inf}
\|v\|_{L^{\infty}(\Om)} \leq  C_1\, 
\bigg(1+\|u\|_{H^1(\Om)}^{B}\bigg)\, 
\bigg(1+\|v\|_{H^1(\Om)}^{B}\bigg)\,,
\end{equation}
where $A,\ B$ depend on $p_1,p_2$ and $N$  are given explicitly by
$$
A:=\frac{1}{(N-1) (p_1+1)\de_0}\,,\qquad
B:=\frac{1}{(N-1) (p_2+1)\de_0}\,.
$$
\end{thm}

\par 
The scalar case of Theorem~\ref{sys:th:apriori:estim}  is considered in \cite{CMP_2024}.   Similar estimates as \eqref{sys:est:u:L:inf}-\eqref{sys:est:v:L:inf} for the scalar case, when the nonlinearity is at the interior, can be found in \cite{Pardo_JFPTA} for the semilinear case,   and \cite{Pardo_RACSAM} for the $p$-Laplacian case, including singular
Carathéodory nonlinearities in both situations. 

\par A uniform $H^1(\Om)$ {\it a priori} bound result for weak solutions of \eqref{sys:u}-\eqref{sys:v} will complement our result, providing a uniform $L^\infty(\Om)$ {\it a priori} bound. Theorem \ref{sys:th:apriori:estim} implies in particular that sequences of weak solutions  uniformly  bounded in  their $H^1(\Om)$-norms, are uniformly bounded in their $L^\infty(\Om)$-norm.
In fact, a question naturally arises: does the coupled system case \eqref{sys:u}-\eqref{sys:v} have a uniform $L^\infty(\Om)$ {\it a priori} bound for any positive weak solution? In \cite[Theorem 3.7]{Bon-Ros_2001}, the authors give an affirmative answer   when $p_1p_2>1$ and 
$p_1,p_2\le 2_*-1$ but not both equal for $N \ge 3$. We conjecture that such a result is true when $p_1, p_2$ satisfy $p_1p_2>1$ and the critical hyperbola condition \eqref{crit:hyp}, which to the best of our understanding,  remains as an open problem.

\par In \cite{Marino-Winkert_systems_2020},
authors consider a general quasi-linear elliptic system with nonlinear boundary conditions. Their work deals with nonlinearities in the interior as well as on the boundary,
depending on both components simultaneously, and establish that the weak solutions are in $\big(L^{\infty}(\Om)\cap L^{\infty}(\p\Om)\big)^2$, see \cite[Thm.~3.1, Thm.~3.2]{Marino-Winkert_systems_2020}. However, their  results do not contain the result of Theorem~\ref{sys:th:reg}. Indeed, their hypotheses (E16) and (E18) applied to \eqref{sys:u}-\eqref{sys:v} and \eqref{f:g:growth} imply that $p_1,\, p_2$ satisfy  $1< p_1, p_2 <2_*-2=\frac{2}{N-2}$, which is valid only for $N < 4$.  Further, they mentioned, see \cite[Remark~3.3]{Marino-Winkert_systems_2020},
 that the natural assumptions will be (E16') and (E18'), which applied to \eqref{sys:u}-\eqref{sys:v} and \eqref{f:g:growth} require that $p_1, p_2 <2_*-1=\frac{N}{N-2}$, 
 which will still be more restrictive than the  critical hyperbola condition \eqref{sys:crit:hyp:2}
 of Theorem~\ref{sys:th:reg}. The $p_1p_2$ regions covered by Theorem~\ref{sys:th:reg} and   the result of \cite{Marino-Winkert_systems_2020} 
 are compared in Figure~\ref{fig:region} for different values of $N\ge 3$.  

\par See  also \cite{Ho_Winkert,  Winkert_2010} and references therein, where regularity results, similar to Theorem~\ref{sys:th:reg} for the scalar case  are investigated for more general operators and nonlinearities. 

\par The  critical hyperbola condition with nonlinearities at the interior first appears independently in \cite{Clement_deFigueiredo_Mitidieri} and \cite{Peletier_vanderVorst}. 
In \cite{Clement_deFigueiredo_Mitidieri} the authors focus on the existence of uniform {\it a priori} bounds below the critical hyperbola, that is, $M_i$ independent of $u$ and $v$. In  \cite{Peletier_vanderVorst}, the authors focus on regions separating existence and nonexistence of radial solutions to elliptic systems.\\

\par  To prove Theorem~\ref{sys:th:reg}, we  first prove a regularity result adapted for nonlinear boundary conditions, Lemma~\ref{sys:Moser},  using a version of the De Giorgi-Nash-Moser  iteration scheme,
see \cite{DeGiorgi, Nash, Moser}. It is an iterative method based on truncation techniques to obtain an $L^\infty(\Om)$ {\it a priori} bounds.
We use Theorem~\ref{sys:th:reg}, elliptic regularity theory combined with Sobolev embeddings and Gagliardo--Nirenberg interpolation theory  \cite{Brezis_2011, Nirenberg} to prove Theorem~\ref{sys:th:apriori:estim}. 

\par Assuming appropriate conditions on the smoothness of $\p\Om$ and on the nonlinearities, we then use the regularity results in \cite[Theorem 2]{Lieberman_1988} to prove further regularity of the weak solution, see Corollary~\ref{sys:lem:reg:2}. 

 \medskip

The paper is organized as follows. Section \ref{sec:prelim} is devoted to some preliminaries. In Section \ref{sec:Moser:sys} we 
state and prove a regularity result 
by using De Giorgi-Nash-Moser iteration type arguments, see Lemma \ref{sys:Moser}. We prove Theorem \ref{sys:th:reg} in Section \ref{sec:proof:sys:reg}, and also state and prove a corollary to Theorem \ref{sys:th:reg},   Corollary \ref{sys:lem:reg:2}. In Section \ref{sec:proof:th:esim:cor}, we prove Theorem \ref{sys:th:apriori:estim}.

\par In Appendix~\ref{app:A}, we compute the $H^1(\Om)$ norm of a test function, and in Appendix-\ref{app:B}, we include the details of some technical estimate important in Section~\ref{sec:Moser:sys}.
In Appendix~\ref{app:reg}, we collect some known regularity results for  general non-homogeneous Neumann and  Dirichlet linear problems, which are scattered in the literature.

\section{Preliminaries}
\label{sec:prelim}
Throughout this paper, we make use of the continuity of the trace operator 
\begin{equation} \label{trace}
 H^1(\Omega) \hookrightarrow L^{2_*}(\p\Om), \quad \text {where}\;\; 2_*:=\frac{2(N-1)}{N-2}
 \end{equation}
and the Sobolev embedding 
 \begin{equation}\label{Sobolev:embed}
 H^1(\Omega)\hookrightarrow L^{2^*}(\Om), \quad \text {where}\;\;2^*:=\frac{2N}{N-2}.
   \end{equation}

A pair of functions $(u,v)$ is a {\it weak solution} to \eqref{sys:u}-\eqref{sys:v}, 
if $u,v\in  H^1(\Om)$ and  $v \in  H^1(\Om)\cap L^{(N-1)(p_2-1)}(\p\Om)$  whenever $p_2> 2_* -1$, and $u,v$ satisfy
\begin{equation}
\label{def:weak:system:u}
\int_{\Omega} \left[\nabla u \nabla \phi + u \phi \right]= \int_{\p\Om} f(x, v)\phi \qquad \forall \phi \in H^1(\Om)\,
\end{equation}
and 
\begin{equation}
\label{def:weak:system:v}
\int_{\Omega} \left[\nabla v \nabla \psi + v \psi \right]= \int_{\p\Om} g(x, u)\psi \qquad \forall \psi \in H^1(\Om)\,.
\end{equation}

The {\it LHS} of \eqref{def:weak:system:u} and \eqref{def:weak:system:v} are well defined since $u, v$ and $\psi$ are in $H^1(\Om)$.
We also note that $u,\, v,\, \phi,\, \psi \in L^{2_*}(\p\Om)$ because of \eqref{trace}.

Now, we will verify that the {\it RHS} of \eqref{def:weak:system:u} and \eqref{def:weak:system:v} are well-defined. Since $1<p_1 \leq p_2$ and \eqref{sys:crit:hyp:2} holds, we see that $p_1 \le 2_*-1$. 
The growth condition  \eqref{f:g:growth}  imply that  $g(\cdot, u(\cdot)) \in L^{\frac{2_*}{p_1}}(\p\Om)$. Using H\"older's inequality, we get
\begin{equation*}
    \int_{\p\Om} |g(x, u(x))\psi(x)|\le \|g(\cdot, u(\cdot))\|_{L^{(2_*)'}(\p\Om)}\|\psi\|_{L^{2_*}(\p\Om)} < +\infty,
    \end{equation*}
since  $p_1 \le 2_*-1$ implies   $(2_*)'=\frac{2_*}{2_*-1}=(\frac{2_*}{p_1})(\frac{p_1}{2_*-1})\le\frac{2_*}{p_1}$.

\par For the integral on the {\it RHS} of \eqref{def:weak:system:u}, 
if $p_2 \le 2_*-1$, we proceed as in \eqref{def:weak:system:v}. 
On the other hand, for $p_2 > 2_*-1$, using $v \in L^{(N-1)(p_2-1)}(\p\Om)$, we get $f(\cdot, v(\cdot)) \in L^{\frac{(N-1)(p_2-1)}{p_2}}(\p\Om)$. Then H\"{o}lder's inequality and the growth condition \eqref{f:g:growth} yield
\begin{equation*}
    \int_{\p\Om} |f(x, v)\phi(x)|\le \|f(\cdot, v(\cdot)) \|_{L^{(2_*)'}(\p\Om)}\|\phi\|_{L^{2_*}(\p\Om)} < +\infty
    \end{equation*}
  since  $p_2 \geq  2_*-1$ implies   $(2_*)'=\frac{2_*}{2_*-1}=\frac{2(N-1)}{N} \leq \frac{(N-1)(p_2-1)}{p_2}$.

Observe from above that $1< p_1 \leq 2_*-1$ but $1 < p_2 < +\infty$.  If $1< p_2 \leq 2_*-1$, $u,v\in  H^1(\Om)$ is sufficient to justify that the  integrals on the {\it RHS} of \eqref{def:weak:system:u} and \eqref{def:weak:system:v} are finite. However, for $p_2> 2_* -1$ we further require  $v \in   L^{(N-1)(p_2-1)}(\p\Om)$  for the integral in \eqref{def:weak:system:v} to be finite.
\vskip.2in

Next, we
consider the linear problem 
\begin{equation}
\label{lbp} 
\left\{ \begin{array}{rcll}
-\Delta  v +v &=&0 & \qquad \mbox{in } \Om\,;  \\
\frac{\p   v}{\p \eta}&=&  \mathit{h} &  \qquad \mbox{on } \p
\Om\,,
\end{array}\right.
\end{equation}
where $h \in L^q(\p\Omega)$ for $q \geq 1$. It is known  that for each $q \geq 1$, \eqref{lbp} has a unique solution in $ W^{1,m}(\Om)$ and 
\begin{equation}
\label{ineq:m:h}
\|v\|_{W^{1,m}(\Om)}\le  C\|\mathit{h}\|_{L^{q}(\p\Om)},\qquad \text{where }\quad 1 \leq  m\le Nq/(N-1)\,,
\end{equation}
see, for instance \cite{Mav-Pardo_2017} for more details, or part {\it (i)} of Proposition \ref{pro:reg:lin:N}.

\section{De Giorgi-Nash-Moser type estimate for systems}
\label{sec:Moser:sys}
Using a De Giorgi-Nash-Moser type iteration, we obtain  the regularity of weak solutions to \eqref{sys:u}-\eqref{sys:v} with atmost linear growth on the boundary. This result plays a crucial role in the proof of Theorem~\ref{sys:th:reg}. 

\begin{lem}
    [De Giorgi-Nash-Moser type estimates]
\label{sys:Moser}
Let $\Om$ be a bounded domain in $\R^N$ 
and $f,\ g:\p\Om\times\R \to  \R$ be  Carathéodory functions such that for non-negative functions $ a_1, a_2 \in  L^{N-1}(\p\Om)$
\begin{equation}
\hspace{-4cm}
\begin{cases}
|f(x, t)| \le a_1(x)(1 + |t|),&\\
|g(x, t)| \le a_2(x)(1 + |t|),&
\end{cases}	
\label{hyp:lin:sys}
\end{equation}
for all $t \in \R$ and for a.e. $x \in  \p\Om$. If $(u, v) $ is a weak solution to the system \eqref{sys:u}-\eqref{sys:v}
then $u,v \in  L^r (\p\Om)$ for any $1 \leq r<\infty$. 
\end{lem}

\begin{proof}
Let $(u,v)$ be a weak solution of \eqref{sys:u}-\eqref{sys:v}. Clearly, since $u, v \in H^1(\Om)$, by the continuity of the trace operator $u, v \in L^r(\p\Om)$ for $1 \leq r \leq 2_*$. 
\par To show $u,v \in L^r(\p\Om)$ for $r > 2_*$, we use De Giorgi-Nash-Moser's iteration scheme (cf. \cite[Lem.~B.3]{Struwe} where the author refers this scheme as the Brezis-Kato-Moser iteration technique) by utilizing the properties of  test functions, constructed by  truncating powers of $u$ and $v$. The proof is carried out in several steps below.\\

\noindent  (I): {\it Choice of test functions.}

\smallskip

\noindent For $s\ge 0,L\ge  0$, define 
\begin{align}
\label{def:vf1}
\vf^{(1)}&=\vf_{s,L}^{(1)} := u \min \{|u|^{2s},L^2 \}   \in  H^1(\Om), \mbox{ and}\\
\vf^{(2)}&=\vf_{s,L}^{(2)} := v \min \{|v|^{2s},L^2 \}   \in  H^1(\Om)\,
\label{def:vf2}
\end{align}
where the gradients  are given by 
\begin{equation}
\label{gradient:phi:1}
\nabla \vf^{(1)} =
\left\{
\begin{array}{lr}
  (2s+1)|u|^{2s}\nabla u;&  |u|^{s} \leq L,\\
\\ L^2 \nabla u;& |u|^s >L 
\end{array}
\right.
\end{equation}
\begin{equation}
\label{gradient:phi:2}
\nabla \vf^{(2)} =
\left\{
\begin{array}{lr}
  (2s+1)|v|^{2s}\nabla v;& \quad |v|^{s} \leq L,\\
\\
L^2 \nabla v;& \quad |v|^s >L\,. 
\end{array}
\right.
\end{equation} Using  $\vf^{(1)}$ and $\vf^{(2)}$ as test functions in 
\eqref{sys:u} and \eqref{sys:v}, respectively, we obtain
\begin{equation}
\label{sys:weak:u}
\int_{\Om}\left[\nabla u \nabla \varphi^{(1)} + u \varphi^{(1)}\right] = \int_{\p\Om}f(x,v)\varphi^{(1)}\,,
\end{equation}
\begin{equation}
\label{sys:weak:v}
\int_{\Om}\left[\nabla v \nabla \varphi^{(2)} + v \varphi^{(2)} \right]= \int_{\p\Om}g(x,u)\varphi^{(2)}\,.
\end{equation}

First, we estimate the {\it LHS} and {\it RHS} of   \eqref{sys:weak:u}.\\
\noindent (i) {\it The LHS of \eqref{sys:weak:u}.}\\
Using \eqref{gradient:phi:1}, {\it LHS} of \eqref{sys:weak:u} is given by
\begin{align}
\label{sys:LHS:u}
& \int_{\Om}\left[\nabla u \nabla \varphi^{(1)} + u \varphi^{(1)}\right] 
= 2s \int_{\Om\cap \{|u|^s \leq L\}}|\nabla u|^2|u|^{2s}
+ \int_{\Om\cap \{|u|^s \leq L\}}|\nabla u|^2|u|^{2s} \nonumber\\
&\qquad\qquad + L^2 \int_{\Om\cap\{|u|^s > L\}}|\nabla u|^2 + \int_{\Om} u^2 \min\{|u|^{2s}, L^2\}\\
&= 2s \int_{\Om\cap \{|u|^s \leq L\}} |\nabla u|^2|u|^{2s}+
\int_{\Om}|\nabla u|^2\min\{|u|^{2s}, L^2\}  + \int_{\Om} u^2 \min\{|u|^{2s}, L^2\}\,. \nonumber
\end{align}

\noindent(ii) {\it  The RHS of \eqref{sys:weak:u}.} \\
We first note that   
\begin{equation}\label{est:min}
t\min\{t^{2s}, L^2\} \le 1 + t^2 \min\{t^{2s}, L^2\} \quad\mbox{ for } t \geq 0\,.
\end{equation}
Indeed, assuming without loss of generality that $L\ge 1$, this follows from the facts that for $0 \leq t \leq 1$, $t\min\{t^{2s}, L^2\} \leq 1$, and for $t>1$,   $t\min\{t^{2s}, L^2\}\le t^2 \min\{t^{2s}, L^2\} $ since $t^2>t$. %\ge\min\{t^{2s}, L^2\}$. 
\par Using the hypothesis \eqref{hyp:lin:sys}, we get  
\begin{align*}
& \int_{\p\Om}f(x, v)\varphi^{(1)} \\
&\le \int_{\p\Om} a_1\big(1+|v|\big)|u|\min \big\{|u|^{2s}, L^2 \big\} \,  \nonumber\\
&\le \frac{1}2\int_{\p\Om} a_1\big(|u|^2+|v|^2\big)\min \big\{|u|^{2s}, L^2 \big\} +\int_{\p\Om} a_1 |u|\min \big\{|u|^{2s}, L^2 \big\}\,  \nonumber\\
&= \int_{\p\Om} a_1 \left( |u|\min \big\{|u|^{2s}, L^2 \big\} +\frac{1}2 |u|^2\min \big\{|u|^{2s},L^2\} \right)\,   \nonumber\\ 
&\qquad +\frac{1}2\int_{\p\Om} a_1|v|^2\min \big\{|u|^{2s}, L^2 \big\} \nonumber\\
&\le \int_{\p\Om} a_1 +\frac{3}2\int_{\p\Om} a_1 |u|^2\min \big\{|u|^{2s}, L^2 \big\}\, +
\frac{1}2\int_{\p\Om} a_1 |v|^2\min \big\{|u|^{2s}, L^2 \big\}\, \,, \nonumber
\end{align*}
where the last inequality follows from the estimate \eqref{est:min}. \\

\noindent(iii) {\it The equation \eqref{sys:weak:u}.}\\
\noindent Combining the {\it LHS} given by (i) and {\it RHS} given by (ii) of \eqref{sys:weak:u}, we get
\begin{align}
& 2s \int_{\Om\cap \{|u|^s \leq L\}} |\nabla u|^2|u|^{2s}
+\int_\Om  |\na u|^2  \min \big\{|u|^{2s}, L^2 \big\}  \, 
\nonumber\\
&\qquad +\int_{\Om} |u|^2\min \big\{|u|^{2s}, L^2 \big\} 
\label{est:weak:soln:sys:u} \\
&\le \int_{\p\Om} a_1 +\frac{3}2\int_{\p\Om} a_1 |u|^2\min \big\{|u|^{2s}, L^2 \big\}\, +
\frac{1}2\int_{\p\Om} a_1 |v|^2\min \big\{|u|^{2s}, L^2 \big\}\, \,.\nonumber
\end{align}

\noindent(iv) {\it The equation \eqref{sys:weak:v}.}\\
Similarly, using \eqref{def:vf2} and \eqref{gradient:phi:2}, we get 
\begin{align*}
 & 2s \int_{|v|^s \leq L} |\nabla v|^2|v|^{2s}
+
\int_\Om  |\na v|^2  \min \big\{|v|^{2s}, L^2 \big\}   \, dx
\nonumber\\
&\qquad +\int_{\Om} |v|^2\min \big\{|v|^{2s}, L^2 \big\} \, \nonumber\\
&\le \int_{\p\Om} a_2 +\frac{3}2\int_{\p\Om} a_2 |v|^2\min \big\{|v|^{2s}, L^2 \big\}\, +
\frac{1}{2}\int_{\p\Om} a_2 |u|^2\min \big\{|v|^{2s}, L^2 \big\}\,.
\nonumber
\end{align*}

\vskip.1in

\noindent (II): {\it $u,v  \in  L^{2(s+1)} (\p\Om)$  imply $u,v \in  L^{2_*(s+1)} (\p\Om)$.} 

\smallskip

\noindent We will apply Lemma \ref{norm:us:H1}, and Remark \ref{rem:A},  to estimate  $\|u \min\{|u|^s, L\}\|_{H^1(\Om)}$
in terms of the {\it LHS} of \eqref{sys:weak:u},
by assuming that $u,v \in  H^1(\Om)\cap L^{2(s+1)}(\p\Om)$. 

\smallskip

To estimate the second integral on the  {\it RHS} of \eqref{est:weak:soln:sys:u}, we split $\p\Om={\p\Om\cap \{a_1\le K\}}\cup {\p\Om\cap \{a_1> K\}}$, for $K \geq 1$, to get
\begin{align*}
&\int_{\p\Om} a_1|u|^2\min \big\{|u|^{2s}, L^2 \big\}  \\  
&\ \le K \int_{\p\Om\cap \{a_1\le K\}}
|u|^2\min \big\{|u|^{2s}, L^2 \big\} 
+\int_{{\p\Om\cap \{a_1> K\}}} a_1|u|^2\min \big\{|u|^{2s}, L^2 \big\} \,.
\end{align*}
Then, since $\min \big\{|u|^{2s}, L^2 \big\}\le |u|^{2s}$ and $u\in L^{2(s+1)}(\partial\Om)$,
\begin{equation*}
\int_{{\p\Om\cap \{a_1\le K\}}} |u|^2\min \big\{|u|^{2s}, L^2 \big\}
\leq  \int_{\p \Om}|u|^{2s+2} \leq \text{const}\,.
\end{equation*}
And, using Hölder inequality and Sobolev embedding $H^1(\Om)\hookrightarrow L^{2_*}(\p\Om)$,
\begin{align}\label{est:normH1u}
&\int_{\p\Om\cap \{a_1> K\}} a_1|u|^2\min \big\{|u|^{2s}, L^2 \big\}    \nonumber \\
&\qquad\le \left(\int_{{\p\Om\cap \{a_1> K\}}}  a_1^{N-1}\right)^\frac1{N-1} \left(\int_{\p\Om}  \big| u\min \big\{|u|^{s}, L \big\}\big|^\frac{2(N-1)}{N-2}\right)^\frac{N-2}{N-1} \nonumber \\
&\qquad =  \e_1 (K) \|u \min\{|u|^s, L\}\|^2_{L^{2_*}(\partial\Om)}\nonumber\\
&\qquad\le c_1\e_1 (K)\|u \min\{|u|^s, L\}\|^2_{H^1(\Om)},
\end{align}
where
\begin{align*}
\e_1 (K) := \left(\int_{{\p\Om\cap \{a_1> K\}}}  a_1^{N-1}\right)^\frac1{N-1}\to 0\qq{as}K \to\infty.
\end{align*}
Next we estimate the third integral on the  {\it RHS} of \eqref{est:weak:soln:sys:u}.  On ${\p\Om\cap \{a_1\le K\}}$ we get 
\begin{align*}
& \int_{{\p\Om\cap \{a_1\le K\}}}a_1|v|^2\min\{|u|^{2s}, L^2\} \leq K \int_{\p\Om}|v|^2|u|^{2s} \nonumber \\
&\qquad \leq K \left(\int_{\p\Om}|v|^{2(s+1)}\right)^{\frac{1}{s+1}}\left(\int_{\p\Om}|u|^{2(s+1)}\right)^{\frac{s}{s+1}}\nonumber\\
&\qquad = K \|v\|_{L^{2(s+1)}(\p\Om)}^2\|u\|_{L^{2(s+1)}(\p\Om)}^{2s}=c_0K < \infty\,.
\end{align*}
On ${\p\Om\cap \{a_1> K\}}$, we get
\begin{align*}
&\int_{{\p\Om\cap \{a_1> K\}}} a_1\,|v|^2\,\min \big\{|u|^{2s}, L^2 \big\} \,
% \right)
\le \left(\int_{{\p\Om\cap \{a_1> K\}}}  a_1^{N-1}\right)^\frac{1}{N-1}\\ 
&\qquad\qquad\times
\left(\int_{\p\Om} |v|^{2\frac{(s+1)(N-1)}{N-2}}
\right)^\frac{N-2}{(s+1)(N-1)} 
\left(\int_{\p\Om}   
|u|^{2s\frac{(s+1)(N-1)}{s(N-2)}}
\right)^\frac{s(N-2)}{(s+1)(N-1)} %\Bigg]
\\
&\qquad\le c_2\e_1 (K)
\big\|\,|v|^{s+1}\big\|_{L^{2_*}(\p\Om)}^\frac2{s+1}
\big\| |u|^{s+1}\big\|_{L^{2_*}(\p\Om)}^\frac{2s}{s+1}\,.
\end{align*}
Combining the previous two, the third integral on the {\it RHS} of \eqref{est:weak:soln:sys:u} is estimated as 
\begin{align}
\label{est:weak:<K}
&\int_{\p\Om} a_1 |v|^2\min \big\{|u|^{2s}, L^2 \big\} \\
&\qquad \le  
c_0K+c_2\e_1 (K)
\big\|\,|v|^{s+1}\big\|_{L^{2_*}(\p\Om)}^\frac2{s+1}
\big\| |u|^{s+1}\big\|_{L^{2_*}(\p\Om)}^\frac{2s}{s+1}\,.\nonumber
\end{align}
Now, we apply Lemma \ref{norm:us:H1}, Remark \ref{rem:A}, \eqref{est:weak:soln:sys:u}, \eqref{est:normH1u}, and \eqref{est:weak:<K}  to get 
\begin{align*}
& \big\|u\min \big\{|u|^{s}, L \big\}\big\|_{H^1(\Om)}^2 \\
&\qquad \le c(1+K)+ cK+c_1\e_1 (K)\|u \min\{|u|^s, L\}\|^2_{H^1(\Om)} \\
&\qquad\qquad+c_2\e_1 (K)
\big\|\,|v|^{s+1}\big\|_{L^{2_*}(\p\Om)}^\frac2{s+1}
\big\| |u|^{s+1}\big\|_{L^{2_*}(\p\Om)}^\frac{2s}{s+1}\,,
\end{align*}
where the constant $c$ depends  on $s$, on the $L^{2 s+2}(\p\Om)$ norm of $u$ and $v$ and on the $L^{N-1}(\p\Om)$ norm of $a$ respectively,  but not on $L$. Using the continuity of the trace operator $H^1(\Om)\hookrightarrow L^{2_*}(\p\Om)$, 
\begin{align*}
&(1-\tilde{c}\e_1 (K))\big\|u\min \big\{|u|^{s}, L \big\}\big\|_{H^1(\Om)}^2 \le c(1+2K)\\
&\qquad\qquad\qquad\qquad +
\tilde{c}\e_1 (K)
\big\|\,|v|^{s+1}\big\|_{H^1(\Om)}^\frac2{s+1}
\big\| |u|^{s+1}\big\|_{H^1(\Om)}^\frac{2s}{s+1}\,.
\end{align*}
Then, since $\e_1 (K) \to 0$ as $K \to \infty$, choosing $K$ such that $\tilde{c}\e_1 (K) = \frac14$, 
we conclude that 
\begin{equation*}
\frac{3}{4}\big\|u\min \big\{|u|^{s}, L \big\}\big\|_{H^1(\Om)}^2
\le c(1 + 2K)
+ \frac{1}{4}\, \big\|\,|v|^{s+1}\big\|_{H^1(\Om)}^\frac2{s+1}
\big\| |u|^{s+1}\big\|_{H^1(\Om)}^\frac{2s}{s+1},
\end{equation*}
independent of $L$. Thus letting $L \to\infty$, we have
\begin{equation}
\label{sys:bdd:u}
\frac{3}{4}\big\||u|^{s+1}\big\|_{H^1(\Om)}^2
\le c(1 + 2K)
+ \frac{1}{4}\, \big\|\,|v|^{s+1}\big\|_{H^1(\Om)}^\frac2{s+1}
\big\| |u|^{s+1}\big\|_{H^1(\Om)}^\frac{2s}{s+1}\,.
\end{equation}
Similarly,  
\begin{equation}
\label{sys:bdd:v}
\frac{3}{4}\big\||v|^{s+1}\big\|_{H^1(\Om)}^2
\le c(1 + 2K)
+ \frac{1}{4}\, \big\|\,|v|^{s+1}\big\|_{H^1(\Om)}^\frac{2s}{s+1}
\big\| |u|^{s+1}\big\|_{H^1(\Om)}^\frac{2}{s+1}\,.
\end{equation}
Multiplying \eqref{sys:bdd:u} and  \eqref{sys:bdd:v}, and rearranging terms we get
\begin{align}
\label{sys:bdd:uv}
& \frac{1}{2}\big\|\,|u|^{s+1}\big\|_{H^1(\Om)}^2\,
\big\|\,|v|^{s+1}\big\|_{H^1(\Om)}^2\leq
C_K \\ &
+ C_K\left(
\big\|\,|u|^{s+1}\big\|_{H^1(\Om)}^{\frac{2s}{s+1}}
\big\|\,|v|^{s+1}\big\|_{H^1(\Om)}^{\frac{2}{s+1}}
+ \big\|\,|u|^{s+1}\big\|_{H^1(\Om)}^{\frac{2}{s+1}}
\big\|\,|v|^{s+1}\big\|_{H^1(\Om)}^{\frac{2s}{s+1}}
\right)\,.\nonumber 
\end{align}
where $C_K:=\max\{(c(1 + 2K))^2,c(1 + 2K)/4\}$. Then  
\begin{equation}
\label{sys:key:0}
\big\|\,|u|^{s+1}\big\|_{H^1(\Om)}\,
\big\|\,|v|^{s+1}\big\|_{H^1(\Om)}
\le \hat{C}_K,
\end{equation}
see the details in Appendix~\ref{app:B}. Combining now \eqref{sys:key:0} with \eqref{sys:bdd:u} or \eqref{sys:bdd:v} 
\begin{equation*}
\big\|\,|u|^{s+1}\big\|_{H^1(\Om)}\le \tilde{C}_K\, \mbox{ and }
\big\|\,|v|^{s+1}\big\|_{H^1(\Om)}\le \tilde{C}_K,
\end{equation*}
hence
$$
|u|^{s+1}, |v|^{s+1}\in  H^1(\Om) \hookrightarrow   L^{2_*} (\p\Om). 
$$
That is, whenever $u,\, v \in  L^{2(s+1)} (\p\Om)\cap H^1(\Om)$ we have that $u,\, v \in  L^{2_*(s+1)} (\p\Om)$.\\

\noindent (III): {\it $u, v \in L^r(\p\Om)$ for $r > 2_*$.} \\
Let $s_1 + 1=\frac{2_*}{2}$. Then, since $u,v \in  L^{2_*} (\p\Om)=L^{2(s_1+1)} (\p\Om)$, part (II) implies that $ u,v \in L^{2_*(s_1+1)} (\p\Om)=L^{\left(\frac{2_*}{2}\right)^2}(\p \Om)$. Now, for a  sequence  defined recursively by $s_i+1=\left(\frac{2_*}{2}\right)^i$ for $i \in \N$, part (II) yields  $ u,v \in L^{2_*(s_i+1)} (\p\Om)=L^{\left(\frac{2_*}{2}\right)^{i+1}}(\p \Om)$.

\par Therefore, $u,v \in L^r(\p\Om)$ for any $1 \leq r< \infty$, 
concluding the proof of the lemma. 
\end{proof}

\smallskip

\section{Proof of Theorem~\ref{sys:th:reg}}
\label{sec:proof:sys:reg}
Let $(u, v)$  be a weak solution to \eqref{sys:u}-\eqref{sys:v}.  
We verify that Lemma~\ref{sys:Moser} holds with
$$
a_1(x):=b_1\frac{(1+|v(x)|^{p_2})}{1+|v(x)|} \; \mbox{and} \; 
a_2(x):=b_2\frac{(1+|u(x)|^{p_1})}{1+|u(x)|} \; \mbox{ for } x \in \p\Om\,,
$$
where $b_i >0$ are as in \eqref{f:g:growth}.
First,  for any $t \geq 0$ and $\sigma>1$, the inequality $\dfrac{1+t^{\sigma}}{1+t} \leq 1+t^{\sigma-1}$ holds. Therefore,  
$a_1(x) \leq b_1(1+|v(x)|^{p_2-1})$ and $a_2(x) \leq b_2(1+|u(x)|^{p_1-1})$ on $\p\Om$. 
Then, using  $v\in L^{(N-1)(p_2-1)}(\p \Om)$,  we have 
\begin{align*}
\int_{\p\Om}a_1^{N-1}
& \leq b_1\int_{\p\Om} (1+|v|^{p_2-1})^{N-1} \\
&\leq  \hat{C} |\p \Om| +  \hat{C}\int_{\p\Om} |v|^{(N-1)(p_2-1)} < \infty\,,
\end{align*}
where we used that for any $q>0$ there exists a constant $C>0$ such that $(1+b)^q \leq C(1 + b^q)$ for any $ b \ge 0$, and the fact that $N\ge 3$. Therefore, $a_1 \in L^{N-1}(\p\Om)$.

\smallskip

\par Similarly, to show $a_2 \in L^{N-1}(\p\Om)$, it is enough to show the boundedness of the integral $ \int_{\p\Om}|u|^{(N-1)(p_1-1)} < \infty$. On one hand, $u \in H^1(\Om)$ implies that $u \in L^{2_*}(\p\Om)=L^{\frac{2(N-1)}{N-2}}(\p\Om)$.
On the other hand,  \eqref{sys:crit:hyp:2} implies that 
$$
\frac{2}{p_1+1} \geq \frac{1}{p_1+1} + \frac{1}{p_2+1} \geq \frac{N-2}{N-1}\,,
$$
which gives
$p_1-1 \leq \frac{2}{N-2}$.
Therefore $(N-1)(p_1-1)\le 2_*$ so that $a_2 \in L^{N-1}(\p\Om)$. By Lemma~\ref{sys:Moser} $u, v \in L^q(\p\Om)$ for any $1\leq q < \infty$. 
\par Then $f(\cdot, v(\cdot)),\, g(\cdot, u(\cdot)) \in  L^q(\p\Om)$  for any  $1 \leq q< \infty$.
Now, fixing $q_0> N-1$, Proposition \ref{pro:reg:lin:N}(i) implies that  $u,\ v\in C^\mu (\overline{\Om})\cap W^{1,m}(\Om)$  for any $0<\mu<1$ and  $1<m<\infty$.
 The estimates \eqref{infty:cont:u} and \eqref{infty:cont:v} then follow from \cite[P.~83]{KufnerJohnFucik_book}).
This  concludes the proof. 
\hfill$\Box$

\bigskip

\subsection{A corollary of Theorem \ref{sys:th:reg}}
Next, we establish further regularity of weak solutions of \eqref{sys:u}-\eqref{sys:v} by assuming appropriate conditions on the smoothness of $\p\Om$ and on the nonlinearities.

\begin{cor}[Improved regularity]
\label{sys:lem:reg:2} 
Let  $f,\ g$ satisfy the hypotheses of Theorem~\ref{sys:th:reg}, and
 $(u, v)$  be  a weak solution of  \eqref{sys:u}-\eqref{sys:v}. Suppose that the boundary $\partial\Om$ is $C^{1,\alpha}$  with $0< \alpha \le 1$.  If there exists a constant $L>0$ such that 
\begin{align}
\big|f(x, s)-f(y, t)\big|&\le L\big(|x-y|^\al+|s-t|^\al\big),\label{f:Holder}\\  
\big|g(x, s)-g(y, t)\big|&\le L\big(|x-y|^\al+|s-t|^\al\big), 
\label{g:Holder}
\end{align}
for all $(x,s), (y,t)\in \p \Om\times [-M_0,M_0]$ with $\alpha\in (0,1]$ and constants $M_0 >0 $ and $L\ge 0$,
then  $u, v \in  C^{1,\be}(\overline\Om)$ for some $\be=\be(\al,N, M_0, L)$.  
\par Moreover, 
if the boundary $\partial\Om\in C^{1,1}$ then  $u, v \in C^{2, \mu}(\Om) \cap C^{1,\be}(\overline\Om)$ for any $\mu\in (0,1)$.
\end{cor}
\begin{proof} It follows from Theorem \ref{sys:th:reg} that if $(u,v)$ is a weak solution of \eqref{sys:u}-\eqref{sys:v}, then $u,v\in C(\overline{\Omega})$. Therefore, there exists  a positive constant $\hat{M}$ such that  for all $(x,t)\in \p\Om \times [-M_0,M_0]$
\begin{equation}\label{f:g:bound}
    |f(x,t)|\le\hat{M} \quad \text{and} \quad |g(x,t)|\le \hat{M}\,.
\end{equation}
Using the fact that $\p\Omega\in C^{1,\alpha}$, and  $f$ and $g$ satisfy \eqref{f:Holder}, \eqref{g:Holder} and \eqref{f:g:bound}, by part (v) of Proposition \eqref{pro:reg:lin:N}, $u,v\in C^{2, \beta}(\Om) \cap C^{1,\be}(\overline\Om)$ for some $\beta \in (0,1)$ depending on $\alpha, L, \hat{M}$ and $N$. This completes the proof. 
\end{proof}   

\section{Proof of Theorem \ref{sys:th:apriori:estim}}
\label{sec:proof:th:esim:cor}

 Let $(u,v)$ be a weak solution of \eqref{sys:u}-\eqref{sys:v}. We establish the explicit estimates in several steps below.\\

\noindent{(i) \it  $W^{1,m}(\Om)$ estimates of $u$ and $v$.}\\
By Theorem~\ref{sys:th:reg}, $u \in L^{p_1q_1}(\p\Om)$ and $v \in L^{p_2q_2}(\p\Om)$ for 
\begin{equation}\label{q1:q2}
  q_i>\max\left\{N-1, \frac{2_*}{p_i}\right\},\qquad i=1,2. 
\end{equation} 
Then, it follows from the growth condition \eqref{f:g:growth}, \eqref{crit:hyp} and the elliptic regularity  \eqref{ineq:m:h} that 
\begin{equation}
\label{sys:ineq:m1:q1}
\|u\|_{W^{1,m_1}(\Om)}\le
C\|f(\cdot, v(\cdot))\|_{L^{q_2}(\p\Om)}\le
C\left(|\p\Om| + \|v\|_{L^{p_2q_2}(\p\Om)}\right)
\,,
\end{equation}
\begin{equation}
\label{sys:ineq:m2:q2}
\|v\|_{W^{1,m_2}(\Om)}\le
C\|g(\cdot, u(\cdot))\|_{L^{q_1}(\p\Om)}\le 
C\left(|\p\Om| + \|u\|_{L^{p_1q_1}(\p\Om)}\right)
\,,
\end{equation}
where $m_1:=\frac{Nq_2}{N-1}>N,\ m_2:=\frac{Nq_1}{N-1}>N$.
Moreover, using the estimate \eqref{infty:cont:u} and \eqref{infty:cont:v} of Theorem~\ref{sys:th:reg},
we have
\begin{eqnarray}
\label{sys:est:norm:u:p1}
\|u^{p_1}\|_{L^{q_1}(\p\Om)}^{q_1}
=  \int_{\p \Om}|u|^{p_1q_1-2_*}|u|^{2_*}
\leq \|u\|_{L^\infty(\Om)}^{p_1q_1-2_*}\|u\|_{L^{2_*}(\p \Om)}^{2_*}\,,
\end{eqnarray}
and 
\begin{eqnarray}
\label{sys:est:norm:v:p2}
\|v^{p_2}\|_{L^{q_2}(\p\Om)}^{q_2}
=  \int_{\p \Om}|v|^{p_2q_2-2_*}|v|^{2_*}
\leq \|v\|_{L^\infty(\Om)}^{p_2q_2-2_*}\|v\|_{L^{2_*}(\p \Om)}^{2_*}\,.
\end{eqnarray}
Using \eqref{sys:est:norm:u:p1} in \eqref{sys:ineq:m1:q1},  we get 
\begin{equation}
\label{sys:ineq:m1:q1:2}
\|u\|_{W^{1,m_1}(\Om)}\le
C\left(1+\|v\|_{L^\infty(\Om)}^{p_2-\frac{2_*}{q_2}}\|v\|_{L^{2_*}(\p \Om)}^{\frac{2_*}{q_2}}\right)\,, 
\end{equation}
for some $C$ depending also on the measure $|\p\Om|$, but independent of $u$ and $v$. Similarly, using \eqref{sys:est:norm:v:p2} in \eqref{sys:ineq:m2:q2}, $v$ satisfies
$$
\|v\|_{W^{1,m_2}(\Om)}\le
C\left(1+\|u\|_{L^\infty(\Om)}^{p_1-\frac{2_*}{q_1}}\|u\|_{L^{2_*}(\p \Om)}^{\frac{2_*}{q_1}}\right)\,. 
$$

\noindent{(ii) \it   Gagliardo-Nirenberg interpolation inequalities.}\\
By Gagliardo-Nirenberg interpolation inequality,  see \cite{Nirenberg}, or \cite[3.C~Ex.~3, p.~313-314]{Brezis_2011}, there exists a constant $C=C(|\Om|, N, m_1)$ such that 
\begin{equation}
\label{sys:ineq:gagliardo:system:1}
\|u\|_{L^{\infty}(\Om)} \leq C\|u\|_{W^{1,m_1}(\Om)}^{\sigma_1} \|u\|_{L^{2^*}(\Om)}^{1-\sigma_1}
\end{equation}
where 
\begin{equation}
\label{sys:id:gagliardo:system:1}
\frac{1-\sigma_1}{2^*}
:=\sigma_1 
\left(
\frac{1}{N} - \frac{1}{m_1}
\right) \mbox{ or equivalently } \frac{1}{\sigma_1}:=1+ \frac{2^*}{N} \left(
1 - \frac{N-1}{q_2}
\right).
\end{equation}
Similarly, $v$ satisfies
\begin{equation*}
\|v\|_{L^{\infty}(\Om)} \leq C\|v\|_{W^{1,m_2}(\Om)}^{\sigma_2} \|v\|_{L^{2^*}(\Om)}^{1-\sigma_2}
\end{equation*}
where 
\begin{equation}
\label{sys:id:gagliardo:system:2}
\frac{1-\sigma_2}{2^*}
:=\sigma_2 
\left(
\frac{1}{N} - \frac{1}{m_2}
\right) \mbox{ or equivalently } \frac{1}{\sigma_2}:=1+ \frac{2^*}{N} \left(
1 - \frac{N-1}{q_1}
\right)\,.
\end{equation}
Substituting \eqref{sys:ineq:m1:q1:2} into \eqref{sys:ineq:gagliardo:system:1}, we get

	\begin{equation}
		\label{sys:est:u:sup}
		\|u\|_{L^{\infty}(\Om)} \leq C\left(1+\|v\|_{L^\infty(\Om)}^{\ (p_2-\frac{2_*}{q_2})}\|v\|_{L^{2_*}(\p \Om)}^{\ \frac{2_*}{q_2}}\right)^{\sigma_1} \|u\|_{L^{2^*}(\Om)}^{1-\sigma_1}\,.
	\end{equation}
	Similarly,
	\begin{equation}
		\label{sys:est:v:sup}
		\|v\|_{L^{\infty}(\Om)} \leq C\left(1+\|u\|_{L^\infty(\Om)}^{(p_1-\frac{2_*}{q_1})}\|u\|_{L^{2_*}(\p \Om)}^{\frac{2_*}{q_1}}\right)^{\sigma_2}\|v\|_{L^{2^*}(\Om)}^{1-\sigma_2}\,.
	\end{equation}
Substituting \eqref{sys:est:v:sup} into \eqref{sys:est:u:sup}, 
we get
\begin{align*}
\|u\|_{L^{\infty}(\Om)} 
&\leq C\Bigg[1+\left(1+\|u\|_{L^\infty(\Om)}^{(p_1-\frac{2_*}{q_1})} \|u\|_{L^{2_*}(\p \Om)}^{\frac{2_*}{q_1}}\right)^{\sigma_2(p_2-\frac{2_*}{q_2})}\\
& \qquad\qquad \times \|v\|_{L^{2^*}(\Om)}^{(1-\sigma_2)
	\ (p_2-\frac{2_*}{q_2})}
\|v\|_{L^{2_*}(\p \Om)}^{\frac{2_*}{q_2}}\Bigg]^{\sigma_1} \|u\|_{L^{2^*}(\Om)}^{1-\sigma_1}\,.
\end{align*}
Next, we use the fact that for any  $a>0$, there exist two positive constants $C_1,\ C_2$ such that 
$$
C_1(1+x^a)\le (1+x)^a \le C_2(1+x^a) \qq{for all} x\ge 0\,.
$$
Then, we deduce
\begin{align}
	\|u\|_{L^{\infty}(\Om)} 
	&\leq C\Bigg[1+\left(1+\|u\|_{L^\infty(\Om)}^{\sigma_2(p_1-\frac{2_*}{q_1})(p_2-\frac{2_*}{q_2})} \|u\|_{L^{2_*}(\p \Om)}^{\sigma_2\frac{2_*}{q_1}(p_2-\frac{2_*}{q_2})}\right)\nonumber\\
	& \qquad\qquad \times \|v\|_{L^{2^*}(\Om)}^{(1-\sigma_2)
		\ (p_2-\frac{2_*}{q_2})}
	\|v\|_{L^{2_*}(\p \Om)}^{\frac{2_*}{q_2}}\Bigg]^{\sigma_1} \|u\|_{L^{2^*}(\Om)}^{1-\sigma_1}\nonumber\\
	&\leq C\Bigg[1+\left(1+\|u\|_{L^\infty(\Om)}^{\sigma_2(p_1-\frac{2_*}{q_1})(p_2-\frac{2_*}{q_2})} \|u\|_{L^{2_*}(\p \Om)}^{\sigma_2\frac{2_*}{q_1}(p_2-\frac{2_*}{q_2})}\right)^{\sigma_1}\nonumber\\
	& \qquad\qquad \times \|v\|_{L^{2^*}(\Om)}^{\sigma_1(1-\sigma_2)
		\ (p_2-\frac{2_*}{q_2})}
	\|v\|_{L^{2_*}(\p \Om)}^{\sigma_1\frac{2_*}{q_2}}\Bigg] \|u\|_{L^{2^*}(\Om)}^{1-\sigma_1}\nonumber\\
	&\leq C\Bigg[1+\left(1+\|u\|_{L^\infty(\Om)}^{\sigma_1\sigma_2(p_1-\frac{2_*}{q_1})(p_2-\frac{2_*}{q_2})} \|u\|_{L^{2_*}(\p \Om)}^{\sigma_1\sigma_2\frac{2_*}{q_1}(p_2-\frac{2_*}{q_2})}\right)^{}\nonumber\\
	& \qquad\qquad \times \|v\|_{L^{2^*}(\Om)}^{\sigma_1(1-\sigma_2)
		\ (p_2-\frac{2_*}{q_2})}
	\|v\|_{L^{2_*}(\p \Om)}^{\sigma_1\frac{2_*}{q_2}}\Bigg] \|u\|_{L^{2^*}(\Om)}^{1-\sigma_1}\,.
\label{asymp:ineq:infty:gagliardo:sys}
\end{align}
To establish the explicit estimates, we consider two cases:  $\|u\|_{L^{\infty}(\Om)} > 1$ and   $\|u\|_{L^{\infty}(\Om)} \le 1$.
\par If $\|u\|_{L^{\infty}(\Om)} > 1$,  dividing both sides of \eqref{asymp:ineq:infty:gagliardo:sys} by $\|u\|_{L^{\infty}(\Om)} ^{\sigma_1\sigma_2(p_1-\frac{2_*}{q_1})(p_2-\frac{2_*}{q_2})}$, 
we get
\begin{align}
\label{sys:est:u:sup:final}
\|u\|_{L^{\infty}(\Om)}^{\eta} &\leq   C\Bigg[1+\left(1+ \|u\|_{L^{2_*}(\p \Om)}^{\sigma_1\sigma_2\frac{2_*}{q_1}(p_2-\frac{2_*}{q_2})}\right)^{}\|v\|_{L^{2^*}(\Om)}^{\sigma_1(1-\sigma_2)
	\ (p_2-\frac{2_*}{q_2})}
\|v\|_{L^{2_*}(\p \Om)}^{\sigma_1\frac{2_*}{q_2}}\Bigg]\nonumber\\
& \qquad\qquad \times  \|u\|_{L^{2^*}(\Om)}^{1-\sigma_1}\,,
\end{align}
and similarly 
\begin{align}
\label{sys:est:v:sup:final}
\|v\|_{L^{\infty}(\Om)}^{\eta} &\leq   C\Bigg[1+\left(1+ \|v\|_{L^{2_*}(\p\Om)}^{\sigma_1\sigma_2\frac{2_*}{q_2}(p_1-\frac{2_*}{q_1})}\right) \|u\|_{L^{2^*}(\Om)}^{(1-\sigma_1)\sigma_2(p_1-\frac{2_*}{q_1})} 
\|u\|_{L^{2_*}(\p \Om)}^{\sigma_2\frac{2_*}{q_1}}\Bigg]\nonumber\\
& \qquad\qquad \times \|v\|_{L^{2^*}(\Om)}^{1-\sigma_2}\,,
\end{align}
where 
\begin{equation}
\label{def:eta}
\eta:=1-\sigma_1\sigma_2\left(p_1-\frac{2_*}{q_1}\right)\left(p_2-\frac{2_*}{q_2}\right)<1\,.
\end{equation}\\
If $\|u\|_{L^{\infty}(\Om)} \le 1$, then
\begin{align}
\label{sys:est:u:sup:final:2}
\|u\|_{L^{\infty}(\Om)} &\leq   C\Bigg[1+\left(1+ \|u\|_{L^{2_*}(\p \Om)}^{\sigma_1\sigma_2\frac{2_*}{q_1}(p_2-\frac{2_*}{q_2})}\right)^{}\|v\|_{L^{2^*}(\Om)}^{\sigma_1(1-\sigma_2)
	\ (p_2-\frac{2_*}{q_2})}
\|v\|_{L^{2_*}(\p \Om)}^{\sigma_1\frac{2_*}{q_2}}\Bigg]\nonumber\\
& \qquad\qquad \times  \|u\|_{L^{2^*}(\Om)}^{1-\sigma_1}\,,
\end{align}
and similarly 
\begin{align}
\label{sys:est:v:sup:final:2}
\|v\|_{L^{\infty}(\Om)} &\leq   C\Bigg[1+\left(1+ \|v\|_{L^{2_*}(\p\Om)}^{\sigma_1\sigma_2\frac{2_*}{q_2}(p_1-\frac{2_*}{q_1})}\right) \|u\|_{L^{2^*}(\Om)}^{(1-\sigma_1)\sigma_2(p_1-\frac{2_*}{q_1})} 
\|u\|_{L^{2_*}(\p \Om)}^{\sigma_2\frac{2_*}{q_1}}\Bigg]\nonumber\\
& \qquad\qquad \times \|v\|_{L^{2^*}(\Om)}^{1-\sigma_2}\,,
\end{align}
We observe that the exponents on the  {\it RHS} of \eqref{sys:est:u:sup:final}-\eqref{sys:est:v:sup:final} and of \eqref{sys:est:u:sup:final:2}-\eqref{sys:est:v:sup:final:2} are positive due to \eqref{q1:q2}.\\
\medskip
\par Next, we will verify that also $\eta > 0$. By definition \eqref{def:eta}, $\eta>0$ if and only if
\begin{align*}
\left(p_1-\frac{2_*}{q_1}\right)\left(p_2-\frac{2_*}{q_2}\right) &<   
\frac1{\sigma_1\sigma_2}\,,
\end{align*}
or, equivalently  
\begin{align}
\label{sys:crit:hyp:s}
p_1p_2 &< p_1 \frac{2_*}{q_2} + p_2 \frac{2_*}{q_1}+ \frac1{\sigma_1\sigma_2}  - \frac{2_*2_*}{q_1q_2}\,. 
\end{align}
Multiplying \eqref{crit:hyp} by $(p_2+1)(p_1+1)$ and rearranging terms,
\begin{equation}\label{sys:crit:hyp:4}
p_1p_2 < \frac{1}{N-2}(p_1+p_2)+\frac{N}{N-2}
\,.
\end{equation}
Identifying coefficients on the {\it RHS} of \eqref{sys:crit:hyp:s} and \eqref{sys:crit:hyp:4},  we choose $q_1= q_2$ such that
\begin{align}
\label{sys:crit:hyp:6}
& \frac{2_*}{q_i}=\frac{1}{N-2},\qq{in other words}q_i=2(N-1), 
\end{align}
(observe that $q_i$ satisfy \eqref{q1:q2}), and using the definitions of $\sigma_1$ and $\sigma_2$ in \eqref{sys:id:gagliardo:system:1} and \eqref{sys:id:gagliardo:system:2}, respectively,
\begin{equation}\label{sys:crit:hyp:11}
\sigma_1=\sigma_2=\frac{N-2}{N-1}\,.
\end{equation}
Consequently, 
\begin{align}
\label{sys:crit:hyp:12}
\frac1{\sigma_1\sigma_2}-\frac{1}{(N-2)^2}&
=\frac{N}{N-2}\,.
\end{align}
With the choices of $q_i$ in \eqref{sys:crit:hyp:6}-\eqref{sys:crit:hyp:11}, and hence of $\sigma_i$ in \eqref{sys:crit:hyp:12}, we get that \eqref{sys:crit:hyp:4} implies \eqref{sys:crit:hyp:s}. Hence, $\eta>0$. This completes part (ii).\\

\noindent{(iii) \it  Estimates in terms of Lebesgue norms.}\\
We notice that for the fixed $q_i=2(N-1)$, we get  $m_i=2N$  and 
$$
\eta=\left(\frac{N-2}{N-1}\right)^2\,\left[\frac{1}{N-2}(p_1+p_2)+\frac{N}{N-2}-
p_2p_1\right]\,. 
$$ 
It can also be verified that
\begin{equation}
\label{de:eta}
\eta= \frac{N-2}{N-1}\, (p_1+1)(p_2+1) \de_0\,,
\end{equation}
where $\delta_0>0$ is given by \eqref{de:0}. 
\medskip

Since $\eta>0$, the inequalities \eqref{sys:est:u:sup:final} and  \eqref{sys:est:v:sup:final} can be expressed as
\begin{align}
\label{sys:est:u:L:inf:11}
\|u\|_{L^{\infty}(\Om)} &\leq   C\Bigg[1+\left(1+ \|u\|_{L^{2_*}(\p \Om)}^{\sigma_1\sigma_2\frac{2_*}{q_1}(p_2-\frac{2_*}{q_2})}\right)^{}\|v\|_{L^{2^*}(\Om)}^{\sigma_1(1-\sigma_2)
	\ (p_2-\frac{2_*}{q_2})}
\|v\|_{L^{2_*}(\p \Om)}^{\sigma_1\frac{2_*}{q_2}}\Bigg]^{\frac{1}{\eta}}\nonumber\\
& \qquad\qquad \times  \|u\|_{L^{2^*}(\Om)}^{\frac{1-\sigma_1}{\eta}}\nonumber\\
&\leq   C\Bigg[1+\left(1+ \|u\|_{L^{2_*}(\p \Om)}^{\sigma_1\sigma_2\frac{2_*}{q_1}(p_2-\frac{2_*}{q_2})}\right)^{\frac{1}{\eta}}\|v\|_{L^{2^*}(\Om)}^{\sigma_1(1-\sigma_2)
	\ (p_2-\frac{2_*}{q_2})\frac{1}{\eta}}
\|v\|_{L^{2_*}(\p \Om)}^{\sigma_1\frac{2_*}{q_2}\frac{1}{\eta}}\Bigg]\nonumber\\
& \qquad\qquad \times  \|u\|_{L^{2^*}(\Om)}^{\frac{1-\sigma_1}{\eta}}\nonumber\\
&\leq   C\Bigg[1+\left(1+ \|u\|_{L^{2_*}(\p \Om)}^{\sigma_1\sigma_2\frac{2_*}{q_1}(p_2-\frac{2_*}{q_2})\frac{1}{\eta}}\right)\|v\|_{L^{2^*}(\Om)}^{\sigma_1(1-\sigma_2)
	\ (p_2-\frac{2_*}{q_2})\frac{1}{\eta}}
\|v\|_{L^{2_*}(\p \Om)}^{\sigma_1\frac{2_*}{q_2}\frac{1}{\eta}}\Bigg]\nonumber\\
& \qquad\qquad \times  \|u\|_{L^{2^*}(\Om)}^{\frac{1-\sigma_1}{\eta}}
\,,
\end{align}
and similarly 
\begin{align}
\label{sys:est:v:L:inf:2}
\|v\|_{L^{\infty}(\Om)} &\leq   C\Bigg[1+\left(1+ \|v\|_{L^{2_*}(\p\Om)}^{\sigma_1\sigma_2\frac{2_*}{q_2}(p_1-\frac{2_*}{q_1})\frac{1}{\eta}}\right) \|u\|_{L^{2^*}(\Om)}^{(1-\sigma_1)\sigma_2(p_1-\frac{2_*}{q_1})\frac{1}{\eta}} 
\|u\|_{L^{2_*}(\p \Om)}^{\sigma_2\frac{2_*}{q_1}\frac{1}{\eta}}\Bigg]\nonumber\\
& \qquad\qquad \times \|v\|_{L^{2^*}(\Om)}^{\frac{1-\sigma_2}{\eta}}\,.
\end{align}
Now substituting $q_i$ and $\sigma_i$ as fixed in \eqref{sys:crit:hyp:6} and \eqref{sys:crit:hyp:11}, respectively, in \eqref{sys:est:u:L:inf:11}, we have 

\begin{align}
\label{sys:est:u:L:inf:2}
\|u\|_{L^{\infty}(\Om)} &\leq  C\Bigg[1+\left(1+ \|u\|_{L^{2_*}(\p \Om)}^{\tilde{A}_1}\right)\|v\|_{L^{2^*}(\Om)}^{\tilde{A}_1}
\|v\|_{L^{2_*}(\p \Om)}^{\tilde{A}_2}\Bigg]  \|u\|_{L^{2^*}(\Om)}^{\tilde{A}_2}
\,,
\end{align}
and similarly 
\begin{align}
\label{sys:est:v:L:inf:3}
\|v\|_{L^{\infty}(\Om)} &\leq   C\Bigg[1+\left(1+ \|v\|_{L^{2_*}(\p\Om)}^{\tilde{B}_1}\right) \|u\|_{L^{2^*}(\Om)}^{\tilde{B}_1} 
\|u\|_{L^{2_*}(\p \Om)}^{\tilde{B}_2}\Bigg] \|v\|_{L^{2^*}(\Om)}^{\tilde{B}_2}\,,
\end{align}
\noindent where 
\begin{equation}\label{def:A:B:tilde}
\begin{array}{l}
\tilde{A}_1
:=\displaystyle\frac{\frac{N-2}{(N-1)^2}(p_2-\frac{1}{N-2})}{\eta},
\qquad    
\tilde{A}_2
:=\displaystyle\frac{1}{(N-1)\eta}\,,
\\[.5cm]
\tilde{B}_1
:=\displaystyle\frac{\frac{N-2}{(N-1)^2}(p_1-\frac1{N-2})}{\eta},     
\qquad    
\tilde{B}_2
:=\displaystyle\frac{1}{(N-1)\eta} 
\,,
\end{array}    
\end{equation}
and $C>0$ is independent of $u$ and $v$,  depending only on  $N$  and $\Om$.  Moreover, the exponents $\tilde{A}_i, \tilde{B}_i$ are positive for all $i=1,2$ since $N \ge 3$.
\\

\noindent{(iv) \it  Completing the proof.}\\
Using Sobolev embedding theorems,   and the fact that 
$$
1+(1+x)y\le (1+x)(1+y) \mbox{ for any } x,y\ge 0
$$ 
we get 
\begin{equation*}
\|u\|_{L^{\infty}(\Om)} \leq   C_0 \bigg(1+\|u\|_{H^1(\Om)}^{\tilde{A}_1+\tilde{A}_2} \bigg)
\bigg(1+\|v\|_{H^1(\Om)}^{\tilde{A}_1+\tilde{A}_2} \bigg)\,,
\end{equation*}
and 
\begin{equation*}
\|v\|_{L^{\infty}(\Om)} \leq  C_1  \bigg(1+\|u\|_{H^1(\Om)}^{\tilde{B}_1+\tilde{B}_2}\bigg)
\bigg(1+\|v\|_{H^1(\Om)}^{\tilde{B}_1+\tilde{B}_2}\bigg)\,,
\end{equation*}
where
\begin{equation*}
\tilde{A}_1+\tilde{A}_2
=\frac{\frac{1}{N-1}+\frac{N-2}{(N-1)^2}(p_2-\frac{1}{N-2})}{\eta}
=\frac{\frac{N-2}{N-1} p_2-\left(\frac{N-2}{N-1}\right)^2(p_2-\frac{1}{N-2})}{\eta},
\end{equation*}
and
\begin{equation*}
\tilde{B}_1+\tilde{B}_2
= \frac{\frac{1}{N-1}+\frac{N-2}{(N-1)^2}(p_1-\frac{1}{N-2})}{\eta}
=\frac{\frac{N-2}{N-1} p_1-\left(\frac{N-2}{N-1}\right)^2(p_1-\frac{1}{N-2})}{\eta}.
\end{equation*}
Observe that $\tilde{A}_1+\tilde{A}_2$, and $\tilde{B}_1+\tilde{B}_2$ depend only on  $N$ and $p_1,\ p_2$. \\

\noindent Rearranging terms, 
\begin{equation}
\label{tilde:Ai:new}
\tilde{A}_1+\tilde{A}_2
=\frac{(N-2)(p_2+1)}{\eta (N-1)^2},
\end{equation}
and
\begin{equation}
\label{tilde:Bi:new}
\tilde{B}_1+\tilde{B}_2
=\frac{(p_1+1)-\frac{1}{N-2}}{\eta (N-1)} .
\end{equation}
Due to  \eqref{de:eta},
the exponents \eqref{tilde:Ai:new} and \eqref{tilde:Bi:new} of \eqref{sys:est:u:L:inf:11}-\eqref{sys:est:v:L:inf:2} can be rewritten only in terms of $\de_0$, $p_1$, $p_2$ and $N$, as 
\begin{equation*}
A=\frac{1}{(N-1) (p_1+1)\de_0},\qquad
B=\frac{1}{(N-1) (p_2+1)\de_0}. 
\end{equation*}
That is,  the estimates \eqref{sys:est:u:L:inf:11} and \eqref{sys:est:v:L:inf:2} are 
\begin{equation*}
\|u\|_{L^{\infty}(\Om)} \leq   C_0\,  
\bigg(1+ \|u\|_{H^1(\Om)}^{A} \bigg) 
\bigg(1+\|v\|_{H^1(\Om)}^{A} \bigg) \,
\end{equation*}
\begin{equation*}
\|v\|_{L^{\infty}(\Om)} \leq  C_1\, 
\bigg(1+\|u\|_{H^1(\Om)}^{B}\bigg) \, 
\bigg(1+\|v\|_{H^1(\Om)}^{B}\bigg) \,,
\end{equation*}
as desired. This completes the proof of Theorem~\ref{sys:th:apriori:estim}.
 \hfill $\Box$

\begin{rem} 
{\rm 
The proof of Theorem \ref{sys:th:apriori:estim}  also provides $L^{\infty}$ estimates of weak solutions $(u,v)$ to \eqref{sys:u}-\eqref{sys:v}, in terms of their $L^{2^*}(\Om)$ and $L^{2_*}(\p\Om)$ norms. Specifically,  from the estimates \eqref{sys:est:u:L:inf:2}  and  \eqref{sys:est:v:L:inf:3}  in the proof of Theorem \ref{sys:th:apriori:estim}, it follows  that
\begin{align*}
\|u\|_{L^{\infty}(\Om)} &\leq  C\Bigg[1+\left(1+ \|u\|_{L^{2_*}(\p \Om)}^{\tilde{A}_1}\right)\|v\|_{L^{2^*}(\Om)}^{\tilde{A}_1}
\|v\|_{L^{2_*}(\p \Om)}^{\tilde{A}_2}\Bigg]  \|u\|_{L^{2^*}(\Om)}^{\tilde{A}_2}
\,,
\end{align*}
and 
\begin{align*}
\|v\|_{L^{\infty}(\Om)} &\leq   C\Bigg[1+\left(1+ \|v\|_{L^{2_*}(\p\Om)}^{\tilde{B}_1}\right) \|u\|_{L^{2^*}(\Om)}^{\tilde{B}_1} 
\|u\|_{L^{2_*}(\p \Om)}^{\tilde{B}_2}\Bigg] \|v\|_{L^{2^*}(\Om)}^{\tilde{B}_2}\,.
\end{align*}
\noindent where  \eqref{def:A:B:tilde} and \eqref{de:eta} give
\begin{equation*} 
\begin{array}{l}
\tilde{A}_1=\displaystyle\frac{p_2-\frac{1}{N-2}}{(N-1)\, (p_1+1)(p_2+1) \de_0},
\quad    
\tilde{A}_2=\displaystyle\frac{1}{(N-2)\, (p_1+1)(p_2+1) \de_0}
\\[.5cm]
\tilde{B}_1=\displaystyle\frac{p_1-\frac1{N-2}}{(N-1)\, (p_1+1)(p_2+1) \de_0},     
\quad    
\tilde{B}_2=\displaystyle\frac{1}{(N-2)\, (p_1+1)(p_2+1) \de_0},
\end{array}    
\end{equation*}
and $C>0$ is independent of $u$ and $v$,  depending only on  $N$  and $\Om$.  Observe that 
$\tilde{A}_1+\tilde{A}_2=A$ and $\tilde{B}_1+\tilde{B}_2=B.$
}
\end{rem}

\section*{Appendices}
\begin{appendices}
 \renewcommand{\thesection}{\Alph{section}}
\numberwithin{equation}{section}

\section{$H^1(\Om)$-norm of  $u \min\{|u|^s, L\}$}
\label{app:A}
\begin{lem}
\label{norm:us:H1}
Assume $u \in  H^1(\Om)\cap L^{2(s+1)}(\p\Om) $ for some $s \geq 0$.
Then 
\begin{align*}
&\|u \min\{|u|^s, L\}\|_{H^1(\Om)}=s(s+2) \int_{\Om\cap \{|u|^s \leq L\}}
|\nabla u|^2\min \big\{|u|^{2s}, L^2 \big\} \\
&\qquad  + \int_{\Om}|\nabla u|^2\min \big\{|u|^{2s}, L^2 \big\}
 +\int_{\Om} |u|^2\min \big\{|u|^{2s}, L^2\}\,.\nonumber
\end{align*}
 \end{lem}
 \begin{proof}
Since 
\begin{equation*}
\nabla \min\{|u|^s, L\}
=\left\{
\begin{array}{lll}
& s\, u |u|^{s-2}\nabla u\,; \quad  &|u|^s \leq L\\
& 0; \quad  & |u|^s > L,
\end{array}
\right.
\end{equation*}
and as in \eqref{gradient:phi:1}, we have 
\begin{equation}
\label{gradient:u-times-min}
\nabla ( u\min\{|u|^s, L\})
=\left\{
\begin{array}{lll}
&  (s+1)|u|^s\nabla u\,; \quad  &|u|^s \leq L\\
& L \nabla u; \quad  & |u|^s > L\,.
\end{array}
\right.
\end{equation}

Then, using \eqref{gradient:u-times-min}, we get
\begin{align}
& \|u \min\{|u|^s, L\}\|^2_{H^1(\Om)} \nonumber \\
=&\int_\Om \Big|\na\big(u\min \big\{|u|^{s}, L \big\}\big)\Big|^2
+\int_{\Om} |u|^2\min \big\{|u|^{2s}, L^2 \big\}\nonumber\\ 
=& (s+1)^2 \int_{\Om\cap \{|u|^s \leq L\}}|\nabla u|^2|u|^{2s}
+ L^2\int_{\Om\cap \{|u|^s> L\}}|\nabla u|^2 \nonumber\\ 
&\qquad 
+\int_{\Om} |u|^2\min \big\{|u|^{2s}, L^2\} \nonumber\\ 
=& (s+1)^2 \int_{\Om\cap \{|u|^s \leq L\}}|\nabla u|^2\min \big\{|u|^{2s}, L^2 \big\} 
+ \int_{\Om\cap \{|u|^s> L\}}|\nabla u|^2\min \big\{|u|^{2s}, L^2 \big\}  \nonumber\\ 
&\qquad +\int_{\Om} |u|^2\min \big\{|u|^{2s}, L^2\} \nonumber\\ 
=& s(s+2) \int_{\Om\cap \{|u|^s \leq L\}}
|\nabla u|^2\min \big\{|u|^{2s}, L^2 \big\} + \int_{\Om}|\nabla u|^2\min \big\{|u|^{2s}, L^2 \big\}\nonumber\\ 
&\qquad +\int_{\Om} |u|^2\min \big\{|u|^{2s}, L^2\}\,.\nonumber
\end{align}
\end{proof}
\begin{rem} 
 \label{rem:A}
 Let us observe that 
 \begin{align*}
&\|u \min\{|u|^s, L\}\|^2_{H^1(\Om)} \le  (s+2) \left( 2s\int_{\Om\cap \{|u|^s\le L\}}
|\nabla u|^2\min \big\{|u|^{2s}, L^2 \big\} \right.\nonumber \\
&\qquad\qquad \qquad \left.
+\int_{\Om}
|\nabla u|^2\min \big\{|u|^{2s}, L^2 \big\} 
+\int_{\Om} |u|^2\min \big\{|u|^{2s}, L^2\} \big\}
\right) 
\nonumber\\
&\qquad\qquad  = (s+2) \left( 
\int_{\Om}\nabla u \nabla \varphi^{(1)} + u \varphi^{(1)}
\right)  
\end{align*} 
thanks to \eqref{sys:LHS:u}, and where $\varphi^{(1)}$ is defined in \eqref{def:vf1}.
 \end{rem}

\section{The estimate (3.16)}
\label{app:B}

We establish the estimate  \eqref{sys:key:0} by contradiction. Suppose there exists a sequence of weak solutions $(u_n, v_n)$ of \eqref{sys:u}-\eqref{sys:v} such that 
\begin{equation}
\label{sys:key:0:n}
\big\|\,|u_n|^{s+1}\big\|_{H^1(\Om)}\,
\big\|\,|v_n|^{s+1}\big\|_{H^1(\Om)}
\to +\infty \quad \mbox{as} \quad n \to \infty\,.
\end{equation}
Setting $x_n:=\big\|\,|u_n|^{s+1}\big\|_{H^1(\Om)}$ and 
$y_n:=\big\|\,|u_n|^{s+1}\big\|_{H^1(\Om)}$ \eqref{sys:key:0:n} simplifies to $x_ny_n \to \infty$, where \eqref{sys:bdd:uv} simplifies to 
\begin{equation}
\label{sys:bdd:uv:n}
\frac{1}{2}x_n^2y_n^2 \leq C_k + \tilde{C}_k 
\left(
x_n^{\frac{2}{s+1}}y_n^{\frac{2s}{s+1}} + 
x_n^{\frac{2s}{s+1}}y_n^{\frac{2}{s+1}}
\right)\,.
\end{equation}
Now we consider three cases:
\medskip

\noindent {\bf Case (I):} Suppose $x_n \to \infty$ and $y_n \to \infty$, and define $t_n:=\frac{x_n}{y_n}$.\\
We consider three sub-cases again.\\
\noindent {\bf Case (I.i):} $t_n\to 0$ as $n \to \infty$.\\
Then it follows from \eqref{sys:bdd:uv:n} that 
\begin{equation*}
\frac{1}{2} \leq \frac{C_k}{x_n^2y_n^2} + \tilde{C}_k 
\left(
\frac{t_n^{\frac{2}{s+1}}}{y_n^2} + 
\frac{t_n^{\frac{2s}{s+1}}}{{y_n^2}}
\right) \to 0 \quad  \mbox{as} \quad n \to \infty\,.
\end{equation*}
a contradiction.\\
\noindent {\bf Case (I.ii):} $t_n\to +\infty$ as $n \to \infty$.\\
Then it follows from \eqref{sys:bdd:uv:n} that 
\begin{equation*}
\frac{1}{2} \leq \frac{C_k}{x_n^2y_n^2} + \tilde{C}_k 
\left(
\frac{1}{t_n^{\frac{2s}{s+1}}y_n^2} + 
\frac{1}{t_n^{\frac{2}{s+1}}y_n^2}
\right) \to 0 \quad \mbox{as} \quad n \to \infty\,.
\end{equation*}
a contradiction.\\
\noindent {\bf Case (I.iii):} $t_n \to a \neq 0$ as $n \to \infty$ for some $a \in \R$.\\
In this case, the contradiction follows from any of the two previous cases.
\noindent {\bf Case II:} Suppose (a) $x_n \leq \text{const.}$ and  $y_n \to \infty$ or  (b) $x_n \to \infty$ and $y_n \leq \text{const.}$ as $n \to \infty$. Observe that  \eqref{sys:bdd:u} and \eqref{sys:bdd:v} can be expressed as 
\begin{align*}
(a')\quad     1 \leq \frac{C_0}{x_n^2} + C_1
\frac{y_n^{\frac{2}{s+1}}}{x_n^{\frac{2}{s+1}}}
\quad\mbox{ and }\quad (b') \quad 1 \leq \frac{C_0}{y_n^2} + C_1
\frac{x_n^{\frac{2}{s+1}}}{y_n^{\frac{2}{s+1}}}\,,
\end{align*}
respectively. Then (a) and (b') or  (b) and (a') lead to  contradictions as $n$ goes to infinity. Hence $x_n y_n$ must be bounded, proving the claim.\\

\section{Regularity of the solutions to some linear problems}
\label{app:reg}

Let us consider the Neumann linear problem with reactions at the interior and on the boundary.
\begin{equation}
\label{lnp} 
\left\{ \begin{array}{rcll}
-\Delta  u +u &=&h_1 & \qquad \mbox{in } \Om\,;  \\
\frac{\p   u}{\p \eta}&=&  h_2 &  \qquad \mbox{on } \p
\Om\,,
\end{array}\right.
\end{equation} 
where $\Om$ is a bounded domain in $\mathbb{R}^N$ with $N \ge 2$. Then the following result holds.
\begin{pro}[Regularity of the solutions to the linear problem \eqref{lnp}]\qquad
\label{pro:reg:lin:N}
There exist    a positive constant $C>0$ independent of $u,h_1,h_2$ such that the following holds:  
\begin{enumerate}
    \item[\rm (i)]  If $\p\Om\in C^{0,1}$,  $h_1\in L^{q}(\Om)$ and $h_2\in L^{q'}(\p\Om)$ with $q,\,q'\ge 1$, then there exists a unique $u\in W^{1,p}(\Om)$  and 
    $$
    \|u\|_{W^{1,p}(\Om)}\le C\left(\|h_1\|_{L^{q}(\Om)}+\|h_2\|_{L^{q'}(\p\Om)}\right), 
    $$ 
    where $p=\min\{\frac{qN}{N-q}, \frac{q'N}{N-1} \}$ if $1\le q < N$, or $p=\min\{q,\frac{q'N}{N-1}\}$ if $ q \ge N$.
    Furthermore,  if  $q> N/2$ and $q'>N-1$, then 
    $$
    \|u\|_{C^{\alpha}(\overline{\Om})}\le C\left(\|h_1\|_{L^{q}(\Om)}+\|h_2\|_{L^{q'}(\p\Om)}\right),
    $$
    where $\al=1-N/p$, ($p>N$).
    {\color{red}$\hookrightarrow$ ????????????????????}

    \item[\rm (ii)]  If $\p\Om\in C^{1,1}$,  $h_1\in C^{\alpha}(\Om)\cap L^{q}(\Om)$ and $h_2\in L^{q'}(\p\Om)$ with $q> N/2$ and $q'>N-1$, then there exists a unique $u\in C^{\alpha}(\overline{\Om})\cap C^{2, \alpha}(\Om)$. 
    \item[\rm (iii)] If $\p\Om\in C^{2,\alpha}$, $h_1\in C^{\alpha}(\overline{\Om})$ and $h_2\in C^{1,\alpha}(\p\Om)$ with $\alpha\in (0,1)$, then there exists a unique solution $u\in C^{2,\alpha}(\overline{\Om})$ and 
    $$
    \|u\|_{C^{2, \alpha}(\Omb)}\le C\left(\|h_1\|_{C^{\alpha}(\overline\Om)}+\|h_2\|_{C^{1,\alpha}(\p\Om)}\right), 
    $$
    where $C$ is a positive constant independent of $u,h_1,h_2$. 
    \item[\rm (iv)] If $\p\Om\in C^{2}$, $h_1\in L^p(\Om)$ and $h_2\in W^{1-\frac{1}{p},p}(\p\Om)$ then $u
    \in W^{2,p}(\Om) $ and 
    $$
    \|u\|_{W^{2,p}(\Om)}\le C\left(\|h_1\|_{L^p(\Om)}+\|h_2\|_{ W^{1-\frac{1}{p},p}(\p\Om)}\right),
    $$ 
    where $C$ is a positive constant independent of $u,h_1,h_2$. 

    \item[\rm (v)] If $\p\Om\in C^{1,\alpha}$ with $\alpha\in (0,1]$,  $h_1\in C^{\alpha}(\Om)$ and $h_2\in C^{\alpha}(\p\Om)\cap L^{\infty}({\p\Om})$ then if $u$ is a bounded weak solution of \eqref{lnp} then $u\in C^{1, \beta}(\overline\Om)\cap C^{2, \beta}(\Om)$, where $\beta$ depends on $\alpha$ and $N$. 
\end{enumerate}
\end{pro}

\begin{proof}
\begin{itemize}
    \item[] 
Proof of part (i): It follows from \cite[Ch.3~Sec.~6]{Lad-Ura_1968} or \cite[Lem.~2.2]{Mav-Pardo_2017} that there exists a unique $u\in W^{1,p}(\Om)$ solving \eqref{lnp}. Now if $p>N$,  using the Sobolev embedding theorem, one has $u\in C^{\alpha}(\overline{\Om})$. Then by applying  \cite[Thm.~6.13]{Gil-Trud_2001} for the corresponding nonhomogeneous Dirichlet problem (see part (i) of Proposition \ref{pro:reg:lin:D}), we have  that $u\in C(\overline{\Om})\cap C^{1, \alpha}(\Om)$, see also \cite{Mav-Pardo_2017}.
 \item[] 
Proof of part (ii): From part (i) we have that $u\in C^{\alpha}(\overline{\Om})$. Since $\p\Omega \in C^{1,1}$, $\Om$ satisfies the exterior sphere condition at every point on the boundary and using the fact that $h_1\in C^{\alpha}(\Om)$, reasoning as above it follows from  \cite[Thm.~6.13]{Gil-Trud_2001} that $u\in C^{\alpha}(\overline{\Om})\cap C^{2,\alpha}({\Om})$.
 \item[] 
Proof of part (iii): See \cite[Page 55]{Ama_1976}  or \cite[Chap.3 Sec.~3]{Lad-Ura_1968}.
\item[] Proof of part (iv): See \cite[Page 55]{Ama_1976} or  \cite[Chap.3 Sec.~9]{Lad-Ura_1968}.
\item[]Proof of part (v): By \cite[Thm.~2]{Lieberman_1988}, one has $u \in C^{1,\beta}(\Omb)$. Then using the bootstrap for the differential equation in $\Om$, we get the desired regularity in $\Om$.

\end{itemize}
\end{proof}

Let us consider the Dirichlet linear problem with reactions at the interior and on the boundary
\begin{equation}
\label{lnpD} 
\left\{ \begin{array}{rcll}
-\Delta  u +u &=&h_1 & \qquad \mbox{in } \Om\,;  \\
u&=&  h_2 &  \qquad \mbox{on } \p
\Om\,,
\end{array}\right.
\end{equation} 
where $\Om$ is a bounded domain in $\mathbb{R}^N$ with $N \ge 2$.

\begin{pro}[Regularity of the solutions to the linear problem \eqref{lnpD}]\qquad
\label{pro:reg:lin:D}
There exist    a positive constant $C>0$ independent of $u,h_1,h_2$ such that the following holds:  
\begin{enumerate}
    \item[\rm (i)] If $\p\Om\in C^{1,1}$ and $h_1\in C^{\alpha}(\Om)$ and $h_2\in C(\p\Om)$ then there exists a unique $u\in C(\overline{\Om})\cap C^{2, \alpha}(\Om)$.
    \item[\rm (ii)] If $\p\Om\in C^{2,\alpha}$ and $h_1\in C^{\alpha}(\overline{\Om})$ and $h_2\in C^{2,\alpha}(\p\Om)$ then there exists a unique solution $u\in C^{2,\alpha}(\overline{\Om})$ and 
    $$\|u\|_{C^{2, \alpha}(\partial\Om)}\le C(\|h_1\|_{C^{\alpha}}(\overline\Om)+\|h_2\|_{C^{2,\alpha}(\p\Om)}). 
    $$

    \item[\rm (iii)] If $\p\Om\in C^{2}$  and $h_1\in L^p(\Om)$ and $h_2\in W^{1-\frac{1}{p},p}(\p\Om)$ then $u
    \in W^{2,p}(\Om) $ and 
    $$
    \|u\|_{W^{2,p}(\Om)}\le C\left(\|h_1\|_{L^p(\Om)}+\|h_2\|_{ W^{2-\frac{1}{p}}_p(\p\Om)}\right),
    $$

    \item[\rm (iv)] If $\p\Om\in C^{1,\alpha}$ with $\alpha\in (0,1]$, and $h_1\in C^{\alpha}(\Om)$ and $h_2\in C^{1,\alpha}(\p\Om)$ then if $u$ is a bounded weak solution of \eqref{lnp} then $u\in C^{1, \beta}(\overline\Om)$, where $\beta$ depends on $\alpha$ and $N$.
\end{enumerate}

\end{pro}

\begin{proof}
\begin{itemize}
\item[] 
Proof of part (i): Since $\p\Omega \in C^{1,1}$, $\Om$ satisfies the exterior sphere condition at every point on the boundary. Then  by applying \cite[Thm.~6.13]{Gil-Trud_2001} to the Dirichlet problem, the result holds.
 \item[] 
Proof of part (ii): See \cite[Thm.~6.15]{Gil-Trud_2001} or \cite[Chap.3 Sec.~1-Sec.~2]{Lad-Ura_1968}.
\item[] Proof of part (iii): See \cite[Page 55]{Ama_1976} or \cite[Chap.3 Sec.~9]{Lad-Ura_1968}.
\item[]Proof of part (iv):  See \cite[Thm.~1]{Lieberman_1988}.

\end{itemize}
\end{proof}
\end{appendices}

\section*{Acknowledgements}
This material is based upon work supported by the National Science Foundation under Grant No. 1440140, while the authors were in residence at the Mathematical Sciences Research Institute in Berkeley, California, during the month of June of 2022.
\par The first author  was  supported by a grant from the Simons Foundation 965180. The second author was supported by James A. Michener Faculty Fellowship.
The third author is supported by grants PID2019-103860GB-I00, and PID2022-137074NB-I00,  MICINN,  Spain, and by UCM, Spain, GR58/08, Grupo 920894.

\end{document}